\input amstex

\documentstyle{amsppt}
\magnification=\magstep1 \baselineskip=18pt \hsize=6truein
\vsize=8truein

\topmatter

\title On the uniqueness of the foliation of spheres of constant mean curvature
in asymptotically flat 3-manifolds
\endtitle

\author Jie Qing and Gang Tian
\endauthor

\leftheadtext{Constant mean curvature surfaces} \rightheadtext{Jie
Qing and Gang Tian}

\address Jie Qing, Dept. of Math., UC,
Santa Cruz, Santa Cruz, CA 95064
\endaddress
\email qing{\@}ucsc.edu \endemail

\address Gang Tian, Dept. of Math., Princeton University,
Princeton, NJ 08544
\endaddress
\email tian{\@}math.princeton.edu \endemail

\abstract In this note we study constant mean curvature surfaces in
asymptotically flat 3-manifolds. We prove that, outside a given
compact subset in an asymptotically flat 3-manifold with positive
mass, stable spheres of given constant mean curvature are unique.
Therefore we are able to conclude that there is a unique foliation
of stable spheres of constant mean curvature in an asymptotically
flat 3-manifold with positive mass.
\endabstract
\footnote{After we posted our first version of this paper, we got
some feedbacks indirectly about our main result compared to previous
results on this topic. In this new version we stated our main result
more precisely and gave more specific references for a clearer
comparison of our main result to the previous ones.}


\endtopmatter

\document

\noindent{\bf 1. Introduction} \vskip 0.1in

In the description of isolated gravitational system in General
relativity a space-like time-slice has the structure of a complete
Riemannian 3-manifold with an asymptotically flat end.  Such
asymptotically flat end is diffeomorphic to $R^3\setminus B_1(0)$
and the metric on it asymptotically approaches the Euclidean metric
near the infinity:
$$
g_{ij} = (1+\frac {2m}r)\delta_{ij} + O(r^{-2}),
$$
where $r$ is the Euclidean distance in $R^3$. The constant $m$ can
be interpreted as the total mass of the isolated system and is
referred to as ADM mass in literature [ADM]. It has also been
established in [B] that with reasonable conditions ADM mass can be
geometrically defined independent of the choices of coordinate
system at infinity.

Often it is better to consider an asymptotically flat end as a
perturbation of the static time-slice of the Schwarzchild
space-time.
Let us start with a precise definition of asymptotically flat
3-manifolds adopted from [HY] for our discussions in this note as
follows:

\proclaim{Definition 1.1} A complete Riemannian 3-manifold $(M, g)$
is said to be an asymptotically flat 3-manifold with mass $m$ if
there is a compact domain $K$ of $M$ such that $M\setminus K$ is
diffeomorphic to $R^3\setminus B_1(0)$ and the metric $g$ in this
coordinate system is given as
$$
g_{ij}(x) = (1+\frac m{2|x|})^4 \delta_{ij} + T_{ij}(x),
$$
for all $x\in R^3\setminus B_1(0)$ with a constant $C$ such that
$$
|\partial^l T_{ij}|(x) \leq C |x|^{-2 - l}, \quad 1\leq l\leq 4,
\tag 1.1
$$
where $\partial$ denotes partial derivatives with respect to the
Euclidean coordinates.
\endproclaim

The existence of a unique foliation of spheres of constant mean
curvature near the end in an asymptotically flat manifold is very
important question. Among many applications, the unique foliation of
spheres of constant mean curvature can be used to construct a
geometrically canonical coordinate system at the infinity of
asymptotically flat end. It can also be used to define a geometric
center of mass for an isolated gravitational system (cf. [HY]). The
existence of such unique foliation of spheres of constant mean
curvature at the asymptotically flat end is also helpful to the
study of Penrose inequality regarding the mass (cf. [Br]). In this
note we show that indeed outside a given compact subset in an
asymptotically flat 3-manifold with positive mass there is a unique
foliation of stable spheres of constant mean curvature. Our main
theorem\footnote{The uniqueness problem addressed here was referred
as the global uniqueness of stable CHC surfaces in [HY]. Their
result on this global uniqueness was stated in Theorem 5.1 in [HY].
They proved that for $q > \frac 12$, if $H$ sufficiently small,
there is a unique stable constant mean curvature surface of mean
curvature $H$ outside $B_{H^{-q}}(0)$.
It has been a long-standing question whether stable constant mean
curvature surfaces are unique outside a fixed compact subset. In the
paragraph after Theorem 5.1 on page 301, Huisken and Yau stated:
``{\it it is an open question whether stable constant mean curvature
surfaces are actually completely unique outside a fixed compact
subset.}" Our main theoerm gives an affirmative answer to this
question.} is

\proclaim{Theorem 1.1} Suppose $(M, g)$ is an asymptotically flat
3-manifold with positive mass. Then there exists a compact domain
$K$ such that stable spheres of given constant mean curvature which
separates the infinity from the compact domain $K$ are unique. Hence
there exists a unique foliation of stable spheres of constant mean
curvature outside the compact domain $K$ in $M$.
\endproclaim

The existence of a foliation of stable spheres of constant mean
curvature near asymptotically flat ends was established by Huisken
and Yau in [HY] (also see [Ye]). Some uniqueness results with
additional assumptions were also proven in [Br] [HY] [Ye]. The major
difficulty of establishing the uniqueness of spheres of given
constant mean curvature is that possible drifting of the spheres of
constant mean curvature presents a hurdle to any useful global a
priori estimates on the curvature. As a matter of fact, the
uniqueness is known if one assumes no drifting (cf. [HY] [Ye]).
Moreover, it was proven in [HY] that, if the drifting was somehow
mild, then the uniqueness holds (cf. Theorem 5.1 on page 301 in
[HY]).

Our main technical contributions can be summarized as follows:
First, as a sharp contrast to the Euclidean space, similar to (5.13)
in [HY], we find the following scale invariant integral which
detects the nonzero mass. Suppose that $N$ is a surface of constant
mean curvature in an asymptotically flat end $(R^3\setminus B_1(0),
g)$ with positive mass $m$. Then
$$
\frac 1{8\pi} \int_N \frac H{|x|}\nu\cdot b d\sigma + \frac
1{4\pi}\int_N \frac { (\nu\cdot x)(\nu\cdot b)}{|x|^3}d\sigma \leq C
m^{-1}r_0^{-1}, \tag 1.2
$$
where $d\sigma$ is induced from the Euclidean metric, $C>0$ is some
constant, $b$ is any vector in $R^3$, $\nu$ is the unit out-going
normal vector of $N$ in $R^3$ with respect to the Euclidean metric,
and
$$
r_0 = \min\{|x|: x\in N \subset R^3\setminus B_1(0)\}, \tag 1.3
$$
provided that
$$
\int_N H^2 d\mu < \infty, \tag 1.4
$$
where $d\mu$ is induced from $g$. Secondly\footnote{In [HY], a
global estimate was sought after (cf. Lemma 5.6 in [HY]), with a
compromise to assume that the inner radius is not smaller than
$H^{-q}$ for $q> \frac 12$. They stated in the paragraph after
Theorem 5.1 (page 21, [HY]) that their assumption on inner radius
``{\it seems to be optimal from a technical point of view}". While
in this paper we do different estimates in three different scales.
Particularly we establish some decay estimate for the intermediate
scales by using an asymptotic analysis developed in [QT].}, we are
able to obtain estimates  (cf. Corollary 4.4 and Corollary 4.5 in
Section 4), which are beyond one individual scale in the blow-down
analysis, via an asymptotic analysis used in an early work of us
[QT]. The blow-down for a surface $N$ of constant mean curvature $H$
with the scale $H$ is defined as,
$$
\tilde N = \{ \frac 12 H x: x\in N \subset R^3\setminus
B_1(0)\}\subset R^3. \tag 1.5
$$
The use of the asymptotic analysis introduced in Section 4 is the
key which allows us to obtain some finer estimates and untangle the
problem that uniform roundness and non-drifting of spheres of
constant mean curvature hinge on each other. More precisely, to
eliminate the possible drifting, one carefully calculates the two
integrals in left-hand side of (1.3) for $\tilde N$,
$$
\frac 1{4\pi} \int_{\tilde N} \frac 1 {|x|}\nu\cdot b d\sigma +
\frac 1{4\pi}\int_{\tilde N} \frac { (\nu\cdot x)(\nu\cdot
b)}{|x|^3}d\sigma \tag 1.6
$$
with some particular choice of $b$. If drifting happened, then the
rescaled surface $\tilde N$ would approach the origin. Then one
evaluates the integrals over three different regions: 1) the part of
$\tilde N$ that is any fixed distance away from the origin; 2) the
part of $\tilde N$ that is near the origin in the scale of $H r_0$;
3) the transition between the above two. We will employ Corollary
4.5 in Section 4 to show the integrals on third region contribute
something negligible. Consequently we are able to prove that the
drifting of stable spheres of constant mean curvature does not
happen at all in an asymptotically flat 3-manifold with positive
mass. Then using the early existence and uniqueness results in [HY]
and [Ye], for instance, Theorem 5.1 in [HY], we may conclude our
main theorem.

It is worthwhile to note that the uniqueness of spheres of given
constant mean curvature outside the horizon in the Schwarzchild
space is an interesting open problem. In his thesis [Br], Bray
proved the coordinate spheres are the unique minimizing surfaces of
given constant mean curvature outside the horizon in Scwarzchild
space, in an attempt to prove the Penrose inequality regarding the
mass by the foliation of constant mean curvature surfaces. Theorem
1.1 in the above particularly implies that the coordinate spheres
are the only stable sphere of constant mean curvature near the
infinity of the Schwarzchild space which separates the infinity from
the horizon.

The paper is organized as follows: In Section 2 we will obtain the
curvature estimates based on the Simons' identity and the smallness
of the integral of the traceless part of the second fundamental
form. In Section 3 we introduce the blow-down analysis in all
scales. In Section 4 we recall the asymptotic analysis from [QT] and
prove a technical proposition. Finally in Section 5 we introduce a
sense of the center of mass and prove our main theorem.

\vskip 0.1in\noindent{\bf 2. Curvature estimates} \vskip 0.1in

First let us recall the Simons' identity [SSY] [Sj] for a
hypersurface $N$ in a Riemannian manifold $(M, g)$ (cf. Lemma 1.3 in
[HY]):
$$
\aligned \Delta h_{ij}  & = \nabla_i\nabla_j H + Hh_{ik}h_{jk} -
|A|^2h_{ij} + H R_{3i3j} - h_{ij} R_{3k3k} \\
& \quad + h_{jk}R_{klil} + h_{ik}R_{kljl} - 2h_{lk}R_{iljk} +
\nabla_j R_{3kik} + \nabla_kR_{3ijk} \endaligned \tag 2.1
$$
where $A =(h_{ij})$ is the second fundamental form for $N$ in $M$,
$H = \text{Tr} A$ is mean curvature, and $R_{ijkl}$ and $\nabla
R_{ijkl}$ are curvature and covariant derivatives of curvature for
$(M, g)$. When $N$ is a constant mean curvature hypersurface, we
rather like to rewrite it as an equation for the traceless part
$\AA$ of $A$, i.e. $\AA = A - \frac 12 H$.
$$
\aligned \Delta \AA_{ij} & = H\AA_{ik}\AA_{jk} - \frac 12 H
|\AA|^2 \delta_{ij} - (|\AA|^2 + \frac 12 H^2)\AA_{ij} \\
& \quad + HR_{3i3j} - \frac 12 H R_{3k3k}\delta_{ij} -
\AA_{ij}R_{3k3k} \\
& \quad + \AA_{jk}R_{klil} + \AA_{ik}R_{kljl}  - 2\AA_{lk}R_{likj}\\
& \quad + \nabla_jR_{3kik} - \nabla_kR_{3ikj}.
\endaligned
\tag 2.2
$$

\proclaim{Lemma 2.1} Suppose that $N$ is a constant mean curvature
surface in an asymptotically flat end $(R^3\setminus B_1(0), g)$.
Then
$$
\aligned -|\AA|\Delta |\AA| & \leq |\AA|^4 + CH|\AA|^3 + C H^2
|\AA|^2 \\
& \quad C |\AA|^2 |x|^{-3} + CH|\AA| |x|^{-3} + C|\AA||x|^{-4}.
\endaligned \tag 2.3
$$
\endproclaim

Note that, in an asymptotically flat end (cf. Definition 1.1 in
Section 1),
$$
|R_{ijkl}| \leq C |x|^{-3}, \quad |\nabla R_{ijkl}| \leq C |x|^{-4}.
\tag 2.4
$$
We refer readers to [HY] for the calculations of curvature of the
Schwarzchild space and asymptotically flat ends.

\proclaim{Lemma 2.2} Suppose that $N$ is a constant mean curvature
surface in an asymptotically flat end $(R^3\setminus B_1(0), g)$.
Then $\int_N H_e^2 d\sigma$ is bounded if and only if $\int_N H^2
d\mu$ is bounded, provided that $r_0$ is sufficiently large.
\endproclaim
\demo{Proof} First one may calculate
$$
H_e = (1+\frac m{2r})^2 H + 2 (1+\frac m{2r})^{-1} \frac m{r^3}
x\cdot \nu + O(r^{-3}), \tag 2.5
$$
where $H_e$ is the mean curvature of $N\subset R^3$ with respect to
the Euclidean metric (cf. Lemma 1.4 in [HY]). Hence
$$
H_e^2 = H^2 + O(r^{-1})H^2 + O(r^{-2})H + O(r^{-3}).
$$
Following Lemma 5.2 in [HY] and the fact that $g$ is quasi-isometric
to the Euclidean metric $|dx|^2$, we have:
$$
\aligned \int_N H_e^2 d\sigma & \leq C \int_N H_e^2 d\mu \leq C
\int_N H^2 d\mu +  C (\int_N H^2d\mu)^\frac 12 (\int_N
r^{-4}d\mu)^\frac 12 + C \int_N r^{-3}d\mu \\
& \leq C \int_N H^2 d\mu\endaligned
$$
and
$$
(1 - C r_0^{-1}) \int_N H^2 d\mu \leq C \int_N H_e^2 d\sigma.
$$
Thus the lemma is proved.
\enddemo

Therefore, following Lemma 1 in [Si], we have

\proclaim{Lemma 2.3} Suppose that $N$ is a constant mean curvature
surface in an asymptotically flat end $(R^3\setminus B_1(0), g)$
with $r_0(N)$ sufficiently large, and that
$$
\int_N H^2 d\mu \leq C.
$$
Then
$$
C_1 H^{-1} \leq \text{diam}(N) \leq C_2 H^{-1}.
$$
\endproclaim

We would like to point out that, if the surface $N$ separates the
infinity from the compact part, i.e. the origin is inside $N\subset
R^3$, then the above lemma implies
$$
C_1 H^{-1} \leq r_1(N) \leq C_2 H^{-1}, \tag 2.6
$$
where the outer radius $r_1(N)$ is defined as
$$
r_1(N) = \max \{|x|: x\in N\subset R^3\setminus B_1(0)\}.
$$
Based on Michael and Simon [MS], one has the following Sobolev
inequality (cf. Lemma 5.6 in [HY]).

\proclaim{Lemma 2.4} Suppose that $N$ is a constant mean curvature
surface in an asymptotically flat end $(R^3\setminus B_1(0), g)$
with $r_0(N)$ sufficiently large, and that
$$
\int_N H^2 d\mu \leq C.
$$
Then
$$
(\int_N f^2 d\mu)^\frac 12 \leq C (\int_N |\nabla f|d\mu + \int_N H
|f| d\mu). \tag 2.7
$$
\endproclaim

Now we are ready to state and prove the main curvature estimates:

\proclaim{Theorem 2.5} Suppose that $(R^3\setminus B_1(0), g)$ is an
asymptotically flat end. Then there exist positive numbers
$\sigma_0$, $\epsilon_0$ and $\delta_0$ such that for any constant
mean curvature surface in the end, which separates the infinity from
the compact part, we have
$$
|\AA|^2 (x) \leq C |x|^{-2} \int_{B_{\delta_0 |x|}(x)} |\AA|^2d\mu +
C |x|^{-4} ,\tag 2.8
$$
provided that
$$
\int_N |\AA|^2d\mu \leq \epsilon_0
$$
and $r_0(N) \geq \sigma_0$. And the corresponding a priori estimates
for all covariant derivatives of curvature also hold consequently.
\endproclaim

\demo{Proof} Recall that
$$
\aligned -|\AA|\Delta |\AA| & \leq |\AA|^4 + CH|\AA|^3 + C H^2
|\AA|^2 \\
& \quad C(|\AA|^2 |x|^{-3} + CH|\AA| |x|^{-3} + C|\AA||x|^{-4}).
\endaligned
$$
Multiply the two sides with $\phi^3$, where $\phi$ is an appropriate
cutoff function of small support, and integrate,
$$
\aligned \int_N - \phi^3|\AA|\Delta |\AA|d\mu & \leq \int_N \phi^3
|\AA|^4d\mu + C \int_N H \phi^3 |\AA|^3 d\mu + C \int_N H^2\phi^3
|\AA|^2d\mu \\ & \quad + C r_0^{-1}\int_N \phi^3 ( |\AA|^2 |x|^{-2}
+ CH|\AA| |x|^{-2} + C|\AA||x|^{-3})d\mu \endaligned
$$
where
$$
\aligned \int_N - \phi^3|\AA|& \Delta |\AA|d\mu = \int_N \nabla
(\phi^3|\AA|)\nabla |\AA|d\mu \\ & = \int_N \phi^2 \nabla (\phi
|\AA|)\nabla |\AA| + \int_N 2\phi^2 |\AA| \nabla \phi \nabla |\AA|
d\mu \\ & = \int_N \phi |\nabla (\phi |\AA|)|^2 d\mu + \int_N \phi
|\AA|\nabla (\phi |\AA|)\nabla \phi d\mu + \int_N 2\phi^2 |\AA|
\nabla \phi \nabla |\AA| d\mu \\
&\geq \frac 34 \int_N \phi |\nabla (\phi |\AA|)|^2 d\mu - C \int_N
\phi |\AA|^2 |\nabla \phi|^2d\mu + \int_N 2\phi^2 |\AA|
\nabla \phi \nabla |\AA| d\mu \\
&\geq \frac 12 \int_N \phi |\nabla (\phi |\AA|)|^2 d\mu - C \int_N
\phi |\AA|^2 |\nabla \phi|^2d\mu
\endaligned
$$
$$
\int_N \phi^3 |\AA|^4d\mu \leq (\int_{\text{supp}(\phi)} |\AA|^2
d\mu)^\frac 12 (\int_N(\phi |\AA|)^6d\mu)^\frac 12
$$
and
$$
\int_N H \phi^3|\AA|^3d\mu \leq (\int_{\text{supp}(\phi)}
H^2d\mu)^\frac 12 (\int_N(\phi |\AA|)^6d\mu)^\frac 12.
$$
For other terms
$$
\int_N H^2\phi^3 |\AA|^2 d\mu \leq C r_M^{-2}
\int_{\text{supp}(\phi)} |\AA|^2d\mu \leq C|x_0|^{-2}
\int_{\text{supp}(\phi)} |\AA|^2 d\mu ,
$$
$$
\int_N \phi^3 |x|^{-2} |\AA|^2 d\mu \leq C|x_0|^{-2}
\int_{\text{supp}(\phi)} |\AA|^2d\mu,
$$
$$
\int_N \phi^3 H |x|^{-2} |\AA| d\mu  \leq C
|x_0|^{-2}(\int_{\text{supp}(\phi)} |\AA|^2 d\mu)^\frac 12,
$$
and
$$
\int_N \phi^2 |x|^{-3} |\AA| d\mu \leq C |x_0|^{-2}
(\int_{\text{supp}(\phi)} |\AA|^2 d\mu)^\frac 12.
$$
Note that, for a given point $x_0$, we may choose the cutoff
function $\phi$ so that it has the suppose of a disk of radius, say,
$\delta_0 |x_0|$ ($\delta_0$ to be determined). Now, combining all
terms, we have
$$
\aligned \int_N \phi|\nabla & (\phi |\AA|)|^2 d\mu \leq 2 (\int_N
|\AA|^2 d\mu)^\frac 12 (\int_N(\phi |\AA|)^6d\mu)^\frac 12 + \\
& \quad  C(\int_{\text{supp}(\phi)} H^2d\mu)^\frac 12 (\int_N(\phi
|\AA|)^6d\mu)^\frac 12 + C|x_0|^{-2} (\int_{\text{supp}(\phi)}
|\AA|^2d\mu)^\frac 12.
\endaligned
$$

Applying the Sobolev inequality with $f=\phi^3 g^3$ where $g=|\AA|$,
we have
$$
\aligned (\int_N & (\phi g)^6 d\mu)^\frac 12  \leq C (3 \int_N (\phi
g)^2 |\nabla (\phi g)| d\mu + \int_N H (\phi g)^3) \\
& \leq C ( \int_N \phi ^3 g^4 d\mu)^\frac 12 (\int_N |\nabla (\phi
g)|^2 \phi d\mu)^\frac 12 + (\int_{\text{supp}(\phi)} H^2
d\mu)^\frac
12 (\int_N (\phi g)^6d\mu)^\frac 12 \\
& \leq C (\int_N g^2 d\mu)^\frac 12 (\int_N (\phi g)^6d\mu)^\frac 12
+ C (\int_{\text{supp}(\phi)}
H^2 d\mu)^\frac 12 (\int_N (\phi g)^6d\mu)^\frac 12 \\
& \quad\quad + C \int_N |\nabla (\phi g)|^2 \phi d\mu .
\endaligned
$$
Thus
$$
(\int_N(\phi |\AA|)^6d\mu)^\frac 12 \leq C
|x_0|^{-2}(\int_{\text{supp}(\phi)} |\AA|^2d\mu)^\frac 12, \tag 2.9
$$
which implies
$$
\int_N (\phi|\AA|)^4 d\mu \leq C|x_0|^{-2} \int_{\text{supp}(\phi)}
|\AA|^2d\mu. \tag 2.10
$$
Note that we have chosen $\delta_0$ small enough so that
$$
C \int_{\text{supp}(\phi)} H^2 d\mu \leq \frac 18
$$
and $\int_N |\AA|^2d\mu \leq \epsilon_0$, where $\epsilon_0$ is
small enough so that
$$
C \int_N |\AA|^2 d\mu \leq \frac 18.
$$

Now we proceed to get the point-wise estimates. First, if we take
$f=u^2$ in the Sobolev inequality, then
$$
\aligned (\int_N u^4 d\mu)^\frac 12 & \leq C (2\int_N |u||\nabla
u|d\mu + \int_N H u^2 d\mu) \\
& \leq C (\int_N u^2d\mu)^\frac 12 (\int_N |\nabla u|^2d\mu)^\frac
12 + C(\int_{\text{supp}(u)}H^2d\mu)^\frac 12 (\int_N u^4
d\mu)^\frac 12.
\endaligned
$$
When $u$ has the support as the cutoff function $\phi$, we have
$$
(\int_N u^4 d\mu)^\frac 12 \leq C (\int_N u^2d\mu)^\frac 12 (\int_N
|\nabla u|^2d\mu)^\frac 12. \tag 2.11
$$
To finish the point-wise estimates we use the following rather
standard estimate:

\proclaim{Lemma 2.6} Suppose that a nonnegative function $v$ in
$L^2$ solves
$$
-\Delta v \leq fv +h \tag 2.12
$$
on $B_{2R}(x_0)$, where
$$
\int_{B_{2R}(x_0)} f^2d\mu \leq C R^{-2}
$$
and $h\in L^2(B_{2R}(x_0))$. And Suppose that
$$
(\int_N u^4 d\mu)^\frac 12 \leq C (\int_N u^2d\mu)^\frac 12 (\int_N
|\nabla u|^2d\mu)^\frac 12
$$
holds for all $u$ with support inside $B_{2R}(x_0)$. Then
$$
\sup_{B_R(x_0)} v \leq C R^{-1} \|v\|_{L^2(B_{2R}(x_0))} + C R
\|h\|_{L^2(B_{2R}(x_0))}.
$$
\endproclaim
\demo{Proof} We will simply use the Moser iteration method. For
convenience, we may rescale so that we are working on $B_2$. The
correct scales would be
$$
v_R (x) = v(Rx), f_R (x) = R^2 f(Rx), \ \text{and} \ h_R = R^2
h(Rx).
$$
Let $k = \|h\|_{L^2(B_2)}$ and $\bar v = v + k$. Multiply the
equation with $\phi^2\bar v^{p-1}$ on the both sides
$$
\int |\nabla (\phi \bar v^\frac p2)|^2  \leq \frac 2p \int f\phi^2
\bar v^p + \int h\phi^2 \bar v^{p-1} + C \int |\nabla \phi|^2 \bar
v^p.
$$
Set $\bar f = \frac hk + f$ we have
$$
\int |\nabla (\phi \bar v^\frac p2)|^2  \leq \frac 2p \int \bar
f\phi^2 \bar v^p + C \int |\nabla \phi|^2 \bar v^p.
$$
Note that $\|\bar f\|_{L^2(B_2)} \leq 1 + \|f\|_{L^2(B_2)}$. By the
assumed Sobolev inequality, we
$$
(\int (\phi \bar v^\frac p2)^4)^\frac 12  \leq C (p\int |\bar
f|\phi^2 \bar v^p)^\frac 12 (\int \phi^2 \bar v^p)^\frac 12 + C(\int
|\nabla \phi|^2 v^p)^\frac 12(\int \phi^2 v^p)^\frac 12.
$$
To handle the first term, we apply H\'{o}lder inequality
$$
\aligned (p\int |\bar f|\phi^2 \bar v^p)^\frac 12 (\int \phi^2 \bar
v^p)^\frac 12 & \leq p^\frac 12 (\int |\bar f|^2)^\frac 14 (\int
\phi^2 \bar v^p)^\frac 12
(\int (\phi \bar v^\frac p2)^4)^\frac 14 \\
&  \leq \frac 1{2C}(\int (\phi \bar v^\frac p2)^4)^\frac 12 + C p
(\int |\bar f|^2)^\frac 12 \int \phi^2 \bar v^p.
\endaligned
$$
Hence
$$
(\int (\phi \bar v^\frac p2)^4)^\frac 12 \leq C(p\|\bar f\|_{L^2}
\int \phi^2 \bar v^p + \int |\nabla \phi|^2 \bar v^p).
$$
Now, for $i=1, 2, \dots, $ let $p = 2^i$ and
$$
\phi = \left\{ \aligned  1 & \quad\forall x \in B_{1+2^{-i}} \\
                         0 & \quad\forall x \notin B_{1+ 2^{-i+1}}
               \endaligned\right.
$$
Then
$$
(\int_{B_{1+2^{-i}}} \bar v^{2^{i+1}})^{2^{-i-1} }\leq C^{2^{-i}}
2^{i2^{-i}} (\int_{B_{1+ 2^{-i+1}}} \bar v^{2^i})^{2^{-i}}.
$$
Thus
$$
\sup_{B_1} v \leq \sup_{B_1} \bar v \leq C^{\sum_{i=1} 2^{-i}}
2^{\sum_{i=1} i2^{-i}}(\int_{B_2} \bar v^2)^\frac 12\leq C
(\|v\|_{L^2(B_2)} + \| h\|_{L^2(B_2)}),
$$
whose scaled version gives the lemma.
\enddemo

To get curvature estimates, we write the equation in such way as
(2.12) that we may apply the above lemma for
$$
f = C(|\AA|^2 + H|\AA| + H^2 + r^{-3}) \ \text{and} \ h = C (Hr^{-3}
+ r^{-4}),
$$
in the light of (2.9) and (2.10).
\enddemo

\vskip 0.1in\noindent{\bf 3. Blow-down analysis}\vskip 0.1in

In order to understand a surface of constant mean curvature $N$ in
an asymptotically flat end $(R^3\setminus B_1(0), g)$, we will need
to blow down the surface in different scales. We first consider, the
blow-down by the scale $H$,
$$
\tilde N = \frac 12 H N = \{ \frac 12 H x: x\in N\}. \tag 3.1
$$
Suppose that there is a sequence of constant mean curvature surfaces
$\{N_i\}$ such that
$$
\lim_{i\to\infty} r_0(N_i) = \infty \ \text{and} \ \lim_{i\to
\infty} \int_{N_i} H^2 d\mu = 16\pi. \tag 3.2
$$
Then, by a similar argument to the proof of Lemma 2.2 in the
previous section, we have
$$
\lim_{i\to\infty} \int_{N_i} H_e^2d\sigma = 16\pi. \tag 3.3
$$
Hence, by the curvature estimates established in the previous
section combining the proof of Theorem 1 in [Si], we have

\proclaim{Lemma 3.1} Suppose that $\{N_i\}$ is a sequence of
constant mean curvature surfaces in a given asymptotically flat end
$(R^3\setminus B_1(0), g)$ and that
$$
\lim_{i\to\infty} r_0(N_i) = \infty \ \text{and} \ \lim_{i\to
\infty} \int_{N_i} H^2 d\mu = 16\pi.
$$
And suppose that $N_i$ separates the infinity from the compact part.
Then, there is a subsequence of $\{\tilde N_i\}$ which converges in
Gromov-Hausdorff distance to a round sphere $S^2_1(a)$ of radius $1$
and centered at $a \in R^3$. Moreover, the convergence is in
$C^\infty$ sense away from the origin.
\endproclaim

>From the above lemma, the difficulty will be to study the
possibility of having the origin lying on the sphere $S^2(a)$, that
is,
$$
\lim_{i\to\infty} r_0(N_i) = \infty, \ \text{and} \
\lim_{i\to\infty} r_0(N_i) H(N_i) = 0. \tag 3.5
$$
Then, in the light of the curvature estimates we obtained in the
previous section, we may use the smaller scale $r_0(N_i)$ to blow
down the surface
$$
\hat N = r_0(N)^{-1} N = \{r_0^{-1}x: x\in N\}. \tag 3.6
$$

\proclaim{Lemma 3.2} Suppose that $\{N_i\}$ is a sequence of
constant mean curvature surfaces in a given asymptotically flat end
$(R^3\setminus B_1(0), g)$ and that
$$
\lim_{i\to\infty} r_0(N_i) = \infty \ \text{and} \ \lim_{i\to
\infty} \int_{N_i} H^2 d\mu = 16\pi.
$$
And suppose that
$$
\lim_{i\to\infty} r_0(N_i)H(N_i) = 0.
$$
Then there is a subsequence of $\{\hat N_i\}$ converges to a 2-plane
at distance $1$ from the origin. Moreover the convergence is in
$C^\infty$ in any compact set of $R^3$.
\endproclaim

As one would expect, the real difficulty is to understand the
behavior of the surfaces $N_i$ in the scales between $r_0(N_i)$ and
$H^{-1}(N_i)$. To start we consider the intermediate scales $r_i$
such that
$$
\lim_{i\to\infty} \frac {r_0(N_i)}{r_i} = 0 \ \text{and} \
\lim_{i\to\infty} r_i H(N_i) =0 \tag 3.7
$$
and blow down the surfaces
$$
\bar {N_i} = r_i^{-1}N = \{r_i^{-1} x: x\in N\}. \tag 3.8
$$

\proclaim{Lemma 3.3} Suppose that $\{N_i\}$ is a sequence of
constant mean curvature surfaces in a given asymptotically flat end
$(R^3\setminus B_1(0), g)$ and that
$$
\lim_{i\to\infty} r_0(N_i) = \infty \ \text{and} \ \lim_{i\to
\infty} \int_{N_i} H^2 d\mu = 16\pi.
$$
And suppose that $\{r_i\}$ are such that
$$
\lim_{i\to\infty} \frac {r_0(N_i)}{r_i} = 0 \ \text{and} \
\lim_{i\to\infty} r_iH(N_i) =0.
$$
Then there is a subsequence of $\{\bar N_i\}$ converges to a 2-plane
at the origin in Gromov-Hausdorff distance. Moreover the convergence
is $C^\infty$ in any compact subset away from the origin.
\endproclaim

\vskip 0.1in \noindent {\bf 4. Asymptotic analysis} \vskip 0.1in

In this section we would like to apply the asymptotic analysis used
in [QT] to obtain some estimate that holds over the whole transition
region between the scales $r_0(N_i)$ and $r_1(N_i)$. But first, let
us revise Proposition 2.1 in [QT] as follows. Let us denote
$$
\| u \|_i^2  = \int_{[(i-1)L, iL]\times S^1} | u |^2 dtd\theta.
$$

\proclaim{Lemma 4.1} Suppose $u \in W^{1,2}(\varSigma, R^k)$
satisfies
$$
\Delta u + A \cdot \nabla u + B \cdot u = h \ \ \text{in
$\varSigma$},
$$
where $\Sigma = [0,3L]\times S^1$. And suppose that $L$ is given and
large. Then there exists a positive number $\delta_0$ such that, if
$$
\|h \|_{L^2 (\varSigma)} \leq \delta_0 \max_{1 \leq i \leq 3} \{ \|
u \|_{i} \}
$$
and
$$
\| A \|_{L^\infty (\varSigma)} \leq \delta_0, \ \ \ \
\|B\|_{L^\infty (\varSigma)} \leq \delta_0 ,
$$
then,

(a). $\| u \|_3 \leq e^{-\frac 12 L} \| u \|_{2}$ implies
     $\| u \|_{2} < e^{-\frac 12 L} \| u \|_{1}$,

(b). $\| u \|_{1} \leq e^{-\frac 12 L} \| u \|_{2}$ implies
     $\| u \|_{2} < e^{-\frac 12 L} \| u \|_{3}$, and

(c). If both $\int_{L \times S^1} u d\theta $ and
        $\int_{2L \times S^1} u d\theta \leq \delta_0
        \max_{1 \leq i \leq 3} \{ \| u \|_{i} \}$, then

either
     $\| u \|_{2} <  e^{-\frac 12 L} \| u \|_{1}$ or
     $\| u \|_{2} < e^{-\frac 12 L} \| u \|_{3}$.
\endproclaim
\demo{Proof} We refer to the proof of Proposition 2.1 in [QT] for
more details. In the proof by contradiction argument one needs to
make sure that the sequence of normalized $u_k$ converges to a
non-zero harmonic function $u$ and the non-zero harmonic function
$u$ violates one of (a)-(c). Interior elliptic estimates give the
strong convergence in the middle section $I_2 = [L,2L]\times S^1$,
which implies that $u$ is not trivially zero. Because, with the
assumption of the proof by contradiction, the middle one is the
largest. Finally $u$ indeed induces a contradiction due to the Fatou
lemma.
\enddemo

We would like to point out that Proposition 2.1 in [QT] is
overstated since it is not correct for $l > 3$. But, in the proof of
Proposition 3.1 in [QT], where Corollary 2.2 is used, one may
replace the shifting cylinder with length $3L$ instead of $5L$. The
proof still works the same, which is, one push to the direction of
growth the cylinder of length $3L$ when Corollary 2.2 in [QT]
applies and it gives the estimates regardless of where one is
stopped applying Corollary 2.2.

Given a surface $N$ in $R^3$, recall from, for example, (8.5) in
[Ka], that
$$
\Delta \nu + |\nabla \nu|^2 \nu = \nabla H_e \tag 4.2
$$
where $\nu$ is the Gauss map from $N \longrightarrow S^2$.  For the
constant mean curvature surfaces in the asymptotically flat end
$(R^3\setminus B_1(0), g)$,
$$
|\nabla H_e|(x) \leq C |x|^{-3}. \tag 4.3
$$
Therefore we consider that the Gauss map of the constant mean
curvature surfaces in the asymptotically flat end $(R^3\setminus
B_1(0), g)$ is an almost harmonic map. Hence we are in a situation
which is very similar to that in [QT]. We will refer readers to [QT]
for rather elementary yet involved analysis since the proof we
present here is some modifications from the proof in [QT]. We will
not carry the indices for the surfaces $N_i$ if it does not cause
any confusion. Set
$$
A_{r_1, r_2} = \{x \in N: r_1 \leq |x| \leq r_2\}. \tag 4.4
$$
$A^0_{r_1, r_2}$ stands for the standard annulus in $R^2$. We are
concerned  with the behavior of $\nu$ on the part $A_{Kr_0(N), s
H^{-1}(N)}$ of $N$ where $K$ will be fixed large and $s$ will be
fixed small. The first difference from [QT] is that, while we had a
fixed domain in [QT], we need the following lemma in order to be in
the position to use Lemma 4.1 in the above.

\proclaim{Lemma 4.2} Suppose that $N$ is a constant mean curvature
surfaces in a given asymptotically flat end $(R^3\setminus B_1(0),
g)$ . Then, for any $\epsilon > 0$ and $L$ fixed, there are
$\epsilon_0, s$ and $K$ such that, if
$$
\int_N |\AA|^2 d\sigma \leq \epsilon_0
$$
and $ Kr_0(N) < r < sH^{-1}(N)$, then $(r^{-1}A_{r, e^{4L}r},
r^{-2}g_e)$ may be represented as $(A^0_{1, e^{4L}}, \bar g)$ where
$$
\|\bar g - |dx|^2 \|_{C^1(A^0_{1, e^{4L}})} \leq \epsilon. \tag 4.5
$$
In other words, in the cylindrical coordinates $(S^1\times [\log r,
4L+ \log r], \bar g_c)$,
$$
\|\bar g_c - (dt^2 + d\theta^2)\|_{C^1(S^1\times [\log r, 4L+\log r
])} \leq \epsilon. \tag 4.6
$$
\endproclaim

This is a consequence of Lemma 3.3 in the previous section. Another
difference from [QT] is that, we are considering maps with tension
fields possibly blowing up at point but no energy concentration,
while in [QT] we were considering almost harmonic maps with
concentration of energy but tension fields uniformly bounded in
$L^2$. In cylindrical coordinates, the tension fields
$$
|\tau (\nu)| = r^2 |\nabla H_e| \leq C r^{-1} = C e^{-t} \tag 4.7
$$
for $t\in [\log(K r_0), \log (sH^{-1}]$. Thus,
$$
\int_{S^1\times [t, t+L]} |\tau (\nu)|^2 dtd\theta \leq C e^{-2t}
\tag 4.8
$$
which decays as needed in the arguments in [QT]. But to get the
growth (or decay) of the energy along the cylinder we first can only
have the estimate (3.8) in [QT]. Then we need to use the Hopf
differential
$$
\Phi = |\partial_t \nu|^2 - |\partial_\theta \nu|^2 - 2\sqrt{-1}
\partial_t \nu\cdot \partial_\theta \nu
$$
and the stationary property, in complex variable
$z=t+\sqrt{-1}\theta$,
$$
\bar \partial \Phi = \partial \nu \cdot \tau (\nu) \tag 4.9
$$
to bound $\int |\partial_t\nu|^2$ by $\int |\partial_\theta \nu|^2$
(cf. [QT] [DT]) as follows:
$$
\int_{S^1\times [t, t+L]} |\partial_t\nu|^2dtd\theta \leq
\int_{S^1\times [t, t+L]}|\Phi|dtd\theta + \int_{S^1\times [t,
t+L]}|\partial_\theta\nu|^2dtd\theta.
$$
By the elliptic estimates (cf. [DT]), we have
$$
\int_{S^1\times [t, t+L]} |\Phi| dtd\theta \leq \int_{N\cap
B_r^c}|\Phi| dtd\theta \leq C (\int_{N\cap B_r^c}|\nabla
\nu|^2d\sigma)^\frac 12 (\int_{N\cap B_r^c} |\tau(\nu)|^2
d\sigma)^\frac 12,
$$
where $N \cap B^c_r$ is the part of $N$ which is outside of $B_r$
and is a disk since $N$ is a sphere topologically. Hence, we have
$$
\int_{S^1\times [t, t+L]} |\Phi| dtd\theta \leq C (\int_{N\cap
B_r^c} |\tau(\nu)|^2 d\sigma)^\frac 12 \leq C e^{-t}. \tag 4.10
$$
Notice that in [QT] we instead used the fact that the tension fields
is uniformly bounded in $L^2$ inside $B_\delta$ (cf. lines between
(3.8) and (3.9) in [QT]). The rest of the proof of Proposition 3.1
in [QT] works with little modifications. Thus we have

\proclaim{Proposition 4.3} Suppose that $\{N_i\}$ is a sequence of
constant mean curvature surfaces in a given asymptotically flat end
$(R^3\setminus B_1(0), g)$ and that
$$
\lim_{i\to\infty} r_0(N_i) = \infty \ \text{and} \ \lim_{i\to
\infty} \int_{N_i} H^2 d\mu = 16\pi.
$$
And suppose that
$$
\lim_{i\to\infty} r_0(N_i)H(N_i) = 0.
$$
Then there exist a large number $K$, a small number $s$ and $i_0$
such that, when $i\geq i_0$,
$$
\max_{I_j} |\nabla \nu| \leq C (\int_{B_{sH^{-1}(N_i)}\cap N_i}
|\nabla \nu|^2d\sigma + r_0^{-1})(e^{-\frac 14 jL} + e^{-\frac
14(n_i - j)L}), \tag 4.11
$$
where
$$
I_j = S^1\times [\log(K r_0(N_i)) + jL, \log (K r_0(N_i)) +(j+1)L]
$$
and
$$
j \in [0, n_i] \ \text{and} \ \log (K r_0(N_i)) + (n_i + 1)L = \log
(s H(N_i)^{-1}).
$$
\endproclaim

This finer analysis improves our understanding of the blow-downs
that we discussed in the previous section. Namely,

\proclaim{Corollary 4.4} Assume the same conditions as Proposition
4.3. Then the limit plane in Lemma 3.2 and the limit plane in Lemma
3.3 are all orthogonal to the vector $a$. In fact, we may choose $s$
small and $i$ large enough so that,
$$
|\nu (x) + a | \leq \epsilon
$$
for all $x \in N_i$ and $|x| \leq  s H^{-1}(N_i)$.
\endproclaim

And we have

\proclaim{Corollary 4.5} Asume the same condition as Proposition
4.3. Let $\nu_i = \nu(p_i)$ for some $p_i \in I_{\frac {n_i}2}$.
Then
$$
\aligned \max_{I_j}& |\nu - \nu_i| \\
\leq C \frac 1{1-e^{-\frac 14 L}}(\int_{B_{sH^{-1}(N_i)}\bigcap N_i}
& |\nabla \nu|^2d\sigma + r_0^{-1}) (e^{-\frac 14 jL} + e^{-\frac 18
n_iL}) \endaligned \tag 4.12
$$
for $j\in [0, \frac 12 n_i]$ and
$$
\aligned \max_{I_j} & |\nu - \nu_i| \\ \leq C \frac 1{1-e^{-\frac 14
L}}(\int_{B_{sH^{-1}(N_i)}\bigcap N_i} & |\nabla \nu|^2d\sigma +
r_0^{-1}) (e^{-\frac 18 n_iL} + e^{-\frac 14 (n_i- j)L})
\endaligned \tag 4.13
$$
for $j\in [\frac 12 n_i, n_i]$.
\endproclaim

The two corollary above will be the key for us to calculate the
integrals in next section to prove our main theorem.

\vskip 0.1in\noindent{\bf 5. Center of mass} \vskip 0.1in

First let us recall that, for any embedded surface $N$ in $R^3$ and
any given vector $b\in R^3$,
$$
\int_N H_e \nu\cdot b d\sigma = 0. \tag 5.1
$$
One may consider this as the first variation of the area of surface
$N_t = N + tb \subset R^3$. On the other hand, if $N$ is a constant
mean curvature surface in the asymptotically flat end $(R^3\setminus
B_1(0), g)$, then
$$
\int_N H \nu\cdot b d\sigma = H \int_N \nu\cdot b d\sigma = 0. \tag
5.2
$$
Since the flux is zero across any surface for a given constant
velocity $b$. Thus, for any constant mean curvature surface in the
asymptotically flat end,
$$
\int_N (H_e-H) \nu\cdot b d\sigma = 0. \tag 5.3
$$
One may calculate and find
$$
H_e - H = m (\frac H{|x|} + 2 \frac {\nu \cdot x}{|x|^3}) +
O(|x|^{-2})H + O(|x|^{-3}). \tag 5.4
$$

\proclaim{Lemma 5.1} Suppose $N$ is a surface of constant mean
curvature in the asymptotically flat end with positive mass $m \neq
0$. And suppose that
$$
\int_N H^2 d\mu < \infty
$$
and $r_0(N)$ is sufficiently large. Then for any given $b$ and for
some $C>0$,
$$
\frac 1{8\pi} \int_N \frac H{|x|}\nu\cdot b d\sigma + \frac
1{4\pi}\int_N \frac { (\nu\cdot x)(\nu\cdot b)}{|x|^3}d\sigma \leq C
m^{-1}r_0^{-1}. \tag 5.5
$$
\endproclaim

\demo{Proof} Simply multiply $\nu\cdot b$ to the both sides of (5.4)
and integrate over the surface $N$, we have,
$$
\frac 1{8\pi} \int_N (H_e-H) \nu\cdot b d\sigma = \frac m{8\pi}
\int_N \frac H{|x|}\nu\cdot b d\sigma + \frac m{4\pi}\int_N \frac {
(\nu\cdot x)(\nu\cdot b)}{|x|^3}d\sigma + O(r_0^{-1}).
$$
Here we used the Lemma 5.2 in [HY]. Then the lemma is proved due to
(5.3).
\enddemo

Now, we are ready to state and prove our main theorem in this note
as follows:

\proclaim{Theorem 5.2} Suppose that $\{N_i\}$ is a sequence of
spheres of constant mean curvature in a given asymptotically flat
end with positive mass $m\neq 0$ and that
$$
\lim_{i\to\infty}r_0(N_i) = \infty \ \text{and} \ \lim_{i\to\infty}
\int_{N_i} H^2 d\sigma = 16\pi.
$$
And suppose that $N_i$ separates the infinity from the compact part.
Then
$$
\lim_{i\to \infty} \frac {r_0(N_i)}{r_1(N_i)} =  1. \tag 5.6
$$
\endproclaim
\demo{Proof} We may apply Lemma 3.1 for the blow-down
$$
\tilde N = \frac 12 H N= \{\frac 12 H x: x\in N\}.
$$
If the surfaces $\tilde N_i$ stay away from the origin, i.e.
$$
0 < C \leq H^{-1}_i r_0(N_i)
$$
for some positive constants $C$, then a subsequence of $\tilde N_i$
converges to a sphere $S^2(a)$ radius $1$ and centered at $a\in R^3$
in $C^\infty$ by the curvature estimates Theorem 2.5 in Section 2.
Also notice that (2.6) implies that the blow-down surfaces $\tilde
N_i$ always stay within a bounded region in $R^3$. On one hand, by
(5.5) in Lemme 5.1, we have
$$
\frac 1{4\pi}\int_{S^2(a)} \frac {\nu\cdot b}{|x|}d\sigma + \frac
1{4\pi} \int_{S^2(a)} \frac {(\nu\cdot x)(\nu\cdot b)}{|x|^3}
d\sigma = 0, \tag 5.7
$$
for any $b$. On the other hand, if take $b = - \frac {a}{|a|}$, then
$$
\frac 1{4\pi}\int_{S^2(a)} \frac {\nu\cdot b}{|x|}d\sigma + \frac
1{4\pi} \int_{S^2(a)} \frac {(\nu\cdot x)(\nu\cdot b)}{|x|^3}
d\sigma  = |a|  \tag 5.8
$$
due to an explicit calculation when the origin is inside. Therefore
$$
a = 0 \ \text{and} \ \lim_{i\to\infty} \frac 12 r_0(N_i) H(N_i) =
\lim_{i\to\infty} \frac {r_0(N_i)}{r_1(N_i)} =  1.
$$
To conclude this is all that can happen we need only to exclude the
case when
$$
\lim_{i\to\infty}H^{-1}_i r_0(N_i) =0. \tag 5.9
$$
Assume otherwise, according to Lemma 3.1, the blow-down sequence
$\tilde N_i$ converges to a unit round sphere $S^2(a)$ centered at
$a\in R^3$ with $|a|= 1$ in Hausdorff topology. We will take $b = -
\frac {a}{|a|}$. From Lemma 5.1, we know
$$
\lim_{i\to\infty}(\int_{\tilde N_i} \frac {\nu \cdot b}{|x|}d\sigma
+ \int_{\tilde N_i} \frac {(\nu\cdot x)(\nu\cdot b)}{|x|^3}d\sigma)
= 0. \tag 5.10
$$
But, we claim, on the other hand,
$$
\lim_{i\to\infty}(\int_{\tilde N_i} \frac {\nu \cdot b}{|x|}d\sigma
+ \int_{\tilde N_i} \frac {(\nu\cdot x)(\nu\cdot b)}{|x|^3}d\sigma)
= 4\pi \tag 5.11
$$
which gives us the contradiction. First, we have from explicit
calculations
$$
\int_{S^2(a)} \frac {\nu\cdot b}{|x|}d\sigma   = \frac 43 \pi, \
\int_{S^2(a)} \frac {(\nu\cdot x)(\nu\cdot b)}{|x|^3}d\sigma  =
\frac 23 \pi. \tag 5.12
$$
The first term in (5.11) is an easy term because the uniform
integrability
$$
\lim_{i\to\infty} \int_{\tilde N_i} \frac {\nu \cdot b}{|x|}d\sigma
= \int_{S^2(a)} \frac {\nu\cdot b}{|x|}d\sigma = \frac 43 \pi. \tag
5.13
$$
To deal with the second term in (5.11), we break up the integral
into three parts. For any fixed small number $s >0$ and large number
$K
> 0$,
$$
\aligned \int_{\tilde N_i} \frac {(\nu\cdot x)(\nu\cdot
b)}{|x|^3}d\sigma  = \int_{\tilde N_i \bigcap B^c_s (0)} & \frac
{(\nu\cdot x)(\nu\cdot b)}{|x|^3}d\sigma  \\
+ \int_{\tilde N_i\bigcap (B_s (0)\setminus B_{KHr_0}(0))} \frac
{(\nu\cdot x)(\nu\cdot b)}{|x|^3} d\sigma & + \int_{\tilde
N_i\bigcap B_{KHr_0}(0)} \frac {(\nu\cdot x)(\nu\cdot
b)}{|x|^3}d\sigma. \endaligned \tag 5.14
$$
Then
$$
\lim_{i\to\infty}\int_{\tilde N_i\bigcap B^c_s(0)} \frac {(\nu\cdot
x )(\nu\cdot b)}{|x|^3}d\sigma = \int_{S^2(a) \bigcap B^c_s} \frac
{(\nu\cdot x)(\nu\cdot b)}{|x|^3}d\sigma \tag 5.15
$$
and
$$
\lim_{i\to\infty}\int_{\tilde N_i\bigcap B_{KH (N_i) r_0(N_i)}(0)}
\frac {(\nu\cdot x )(\nu\cdot b)}{|x|^3}d\sigma = \int_{P \bigcap
B_K(0)} \frac {(\nu\cdot x)(\nu\cdot b)}{|x|^3}d\sigma,
$$
where $P$ is the limit plane in Lemma 3.2. By Corollary 4.4, we know
$$
\int_P \frac {(\nu\cdot x)(\nu\cdot b)}{|x|^3}d\sigma = 2\pi.
$$
due to a simple calculation. Notice that
$$
\int_{\tilde N_i} \frac {\nu\cdot x}{|x|^3}d\sigma  = 4\pi \tag 5.16
$$
for any $i$ and
$$
\int_{S^2(a)} \frac {\nu\cdot x}{|x|^3}d\sigma = 2\pi \tag 5.17
$$
because the origin is on the sphere $S^2(a)$. Since
$$
\lim_{i\to\infty}\int_{\tilde N_i\bigcap B^c_s(0)} \frac {\nu\cdot
x}{|x|^3}d\sigma = \int_{S^2(a)\bigcap B^c_s(0)} \frac {\nu\cdot
x}{|x|^3}d\sigma, \tag 5.18
$$
$$
\lim_{i\to\infty}\int_{\tilde N_i\bigcap B_{KH r_0}(0)}\frac
{\nu\cdot x}{|x|^3}d\sigma = \int_{P \bigcap B_K(0)} \frac {\nu\cdot
x}{|x|^3} d\sigma \tag 5.19
$$
and
$$
\int_P \frac {\nu\cdot x}{|x|^3} d\sigma = 2\pi, \tag 5.20
$$
we know
$$
\lim_{i\to\infty, s\to 0, K \to\infty}\int_{\tilde N_i\bigcap (
B^c_s(0)\setminus B_{KHr_0}(0)} \frac {\nu\cdot x}{|x|^3}d\sigma =
0. \tag 5.21
$$
Now we are ready to handle the difficult term: the integral over the
transition region in (5.14). Our goal is to show that
$$
\lim_{i\to\infty, s\to 0, K \to\infty} \int_{\tilde N_i\bigcap (B_s
(0)\setminus B_{KHr_0}(0))} \frac {(\nu\cdot x)(\nu\cdot b)}{|x|^3}
d\sigma = 0. \tag 5.22
$$
The key point is to use Corollary 4.5 to prove (5.22) from (5.21).
Let $\nu_i$ be chosen as in Corollary 4.5. Then
$$
\aligned \int_{\tilde N_i\bigcap (B_s \setminus B_{KHr_0})}
& \frac {(\nu\cdot x)(\nu\cdot b)}{|x|^3}   =  \\
(\nu_i \cdot b) \int_{\tilde N_i\bigcap (B_s \setminus B_{KHr_0})}
\frac {\nu\cdot x}{|x|^3}  & + \int_{\tilde N_i\bigcap (B_s
\setminus B_{KHr_0})} \frac {(\nu\cdot x)((\nu-\nu_i)\cdot
b)}{|x|^3}.
\endaligned\tag 5.23
$$
Hence we only need to deal with the second term on the right side of
the above (5.23). We are better now to use the cylindrical
coordinates used in Section 4.
$$
\aligned & \int_{\tilde N_i\bigcap (B_s (0)\setminus B_{KHr_0}(0))}
\frac {(\nu\cdot x)((\nu-\nu_i)\cdot b)}{|x|^3} d\sigma  \\ & =
\sum_{j=1}^{n_i} \int_{I_j} \frac {(\nu\cdot x)((\nu-\nu_i)\cdot
b)}{|x|^3}  A(t) d\theta dt \\ & \leq C \sum_{j=1}^{n_i} L
\max_{I_j} |\nu - \nu_i| \\
& = C \sum_{j=1}^{n_i /2}L \max_{I_j} |\nu - \nu_i| + C \sum_{j=
n_i/2 +1}^{n_i /2}L \max_{I_j} |\nu - \nu_i| .\endaligned \tag 5.24
$$
>From (4.12) and (4.13), we have
$$
\aligned & \int_{\tilde N_i\bigcap (B_s (0)\setminus B_{KHr_0}(0))}
\frac {(\nu\cdot x)((\nu-\nu_i)\cdot b)}{|x|^3} d\sigma \\ & \leq
C\eta(\sum_{j=1}^{n_i/2}(e^{-\frac 14 Lj} + e^{-\frac 18 Ln_i}) +
\sum_{j=1}^{n_i /2}(e^{-\frac 18 Ln_i} +
e^{-\frac 14 (n_i-j)L})) \\
& \leq C \eta (n_ie^{-\frac 18 n_iL} + 2),
\endaligned \tag 5.25
$$
where
$$
\eta = \int_{B_{sH^{-1}(N_i)}\cap N_i} |\nabla \nu|^2d\sigma +
r_0^{-1}
$$
and $\eta$ can be arbitrarily small as long as $s\to 0$ and
$r_0(N_i)\to\infty$. Thus (5.22) is proved and the proof of the
theorem is completed.
\enddemo

\proclaim{Corollary 5.3} Suppose $(R^3\setminus B_1(0), g)$ is an
asymptotically flat end with positive mass. Then there exist a large
number $K>0$ and a small number $\epsilon >0$ such that, for any $H
< \epsilon$, there exists a unique stable spheres $N$ of constant
mean curvature $H$ with $N \subset R^3\setminus B_K(0)$ and
separates the infinity from the compact part. Hence there exists a
unique foliation of stable spheres of constant mean curvature near
the infinity.
\endproclaim
\demo{Proof} In the light of Proposition 5.3 in [HY] we know
$$
\lim_{i\to\infty} \int_{N_i} H^2 d\mu  = 16\pi, \tag 5.26
$$
provided each $N_i$ is a stable sphere of constant mean curvature.
Thus Corollary 5.3 follows from Theorem 5.2 in the above, Theorem
4.1 in [HY] and the proof of Theorem 5.1 in [HY].
\enddemo

\vskip 0.1in \noindent {\bf References}:

\roster \vskip0.1in
\item"{[ADM]}" R. Amowitt, S. Deser and C. Misner, Coordinate
invariant and energy expression in general relativity, Phys. Rev.
122 (1966) 997-1006.

\vskip0.1in
\item"{[B]}" R. Barnik, The mass of an asymptotically flat
manifold, Comm. Pure. Appl. Math. 39 (1986) 661-693.

\vskip0.1in
\item"{[Br]}" H. Bray, The Penrose inequality in general
relativity and volume comparison theorem involving scalar cutvature,
Stanford Thesis 1997.

\vskip0.1in
\item"{[MS]}" J. Michael and L. Simon, Sobolev and mean-value
inequalities on general submanifolds of $R^n$, Comm. Pure. Appl.
Math. 26 (1973) 361-379.

\vskip0.1in
\item"{[HY]}" G. Huisken and S. T. Yau, Definition of center of
mass for isolated physical systems and unique foliations by stable
spheres of constant mean curvature, Invent. Math. 124 (1996)
281-311.

\vskip0.1in
\item"{[Ka]}" K. Kenmotsu, ``Surfaces with constant mean
curvature'', Tanslations of Math. Monographs, v. 221, AMS, 2003.

\vskip 0.1in
\item"{[QT]}" Jie Qing
and Gang Tian, Bubbling of the heat flow for harmonic maps from
surfaces, Comm. Pure and Appl. Math., Vol. L (1997), 295-310.

\vskip0.1in
\item"{[SSY]}" R. Schoen, L. Simon and S.T. Yau, Curvature
estimates for minimal hypersurfaces, Acta Math. 134 (1975) 275-288.

\vskip0.1in
\item"{[Si]}" L. Simon, Existence of Willmore surface, Proc.
Centre for Math. Anal. 10 (1985) 187-216.

\vskip0.1in
\item"{[Sj]}" J. Simons, Minimal Varieties in Riemannian manifolds,
Ann. of Math. 88 (1968) 62-105.

\vskip0.1in
\item"{[Ye]}" R. Ye, Foliation by constant mean curvature spheres
on asymptotically flat manifolds, Geom. Anal. and the Calculus of
Variations, 369-383, Intern. Press, Cambridge, MA, 1996.

\endroster
\enddocument